\documentclass{amsart}

\usepackage{amsfonts}
\usepackage{amsmath}
\usepackage{amsthm}
\usepackage{amssymb}
\usepackage[english]{babel}
\usepackage{graphics}
\usepackage{latexsym}
\usepackage{longtable} \usepackage{mathrsfs}
\usepackage{graphicx}
\usepackage{wrapfig} 
\usepackage[pdftex,colorlinks=true]{hyperref}

\newtheorem{theorem}{Theorem}
\newtheorem{corollary}[theorem]{Corollary}
\newtheorem{lemma}[theorem]{Lemma}
\newtheorem{proposition}[theorem]{Proposition}

\newtheorem{claim}[theorem]{Claim}
\newtheorem{example}[theorem]{Example}

\theoremstyle{definition}
\newtheorem{definition}[theorem]{Definition}
\newtheorem{remark}[theorem]{Remark}

\renewcommand{\S}{\mathcal{S}}

\newcommand{\A}{\textbf{A}}

\newcommand{\R}{\mathbb{R}}
\newcommand{\C}{\mathbb{C}}
\newcommand{\N}{\mathbb{N}}
\newcommand{\mS}{\mathbb{S}}

\newcommand{\X}{\textbf{X}}
\newcommand{\Y}{\textbf{Y}}
\newcommand{\Z}{\textbf{Z}}

\newcommand{\noi}{\noindent}
\newcommand{\ms}{\medskip}

\newcommand{\al}{\alpha}
\newcommand{\be}{\beta}
\newcommand{\ga}{\gamma}
\newcommand{\de}{\delta}
\newcommand{\De}{\Delta}
\newcommand{\e}{\varepsilon}
\newcommand{\si}{\sigma}
\newcommand{\la}{\lambda}
\newcommand{\ka}{\kappa}
\newcommand{\Om}{\Omega}

\newcommand{\lharpoonup}{-\!\!\!\!\rightharpoonup}

\newcommand{\larrow}{\longrightarrow}
\newcommand{\ot}{\otimes}

\newcommand{\ri}{\rightarrow}

\newcommand{\p}{\partial}
\newcommand{\sub}{\subseteq}

\newcommand{\by}{\times}

\newcommand{\sgn}{\textrm{sgn}}
\newcommand{\ess}{\textrm{ess}}

\newcommand{\Div}{\textrm{Div}}

\newcommand{\cof}{\textrm{cof}}

\newcommand{\bt}{\begin{theorem}}\newcommand{\et}{\end{theorem}}
\newcommand{\bd}{\begin{definition}}\newcommand{\ed}{\end{definition}}
\newcommand{\bl}{\begin{lemma}}\newcommand{\el}{\end{lemma}}
\newcommand{\beq}{\begin{equation}}\newcommand{\eeq}{\end{equation}}
\newcommand{\bc}{\begin{claim}}\newcommand{\ec}{\end{claim}}
\newcommand{\bex}{\begin{example}}\newcommand{\eex}{\end{example}}
\newcommand{\bcor}{\begin{corollary}}\newcommand{\ecor}{\end{corollary}}
\newcommand{\bp}{\begin{proof}}\newcommand{\ep}{\end{proof}}

\newcommand{\BPL}{\medskip \noindent \textbf{Proof of Lemma} }

\newcommand{\BPP}{\medskip \noindent \textbf{Proof of Proposition} }
\newcommand{\BPT}{\medskip \noindent \textbf{Proof of Theorem} }

\numberwithin{equation}{section}

\begin{document}

\title[Existence and Uniqueness for Fully Nonlinear Systems]{Existence and Uniqueness of Global Strong Solutions to Fully Nonlinear Second Order Elliptic Systems}

\author{Nikos Katzourakis}
\address{Department of Mathematics and Statistics, University of Reading, Whiteknights, PO Box 220, Reading RG6 6AX, Berkshire, UK}
\email{n.katzourakis@reading.ac.uk}


\date{}



\begin{abstract} We consider the problem of existence and uniqueness of strong a.e.\ solutions $u: \R^n \longrightarrow \mathbb{R}^N$ to the fully nonlinear PDE system
\[\label{1} \tag{1}
F(\cdot,D^2u ) \,=\,  f, \ \ \text{ a.e.\ on }\mathbb{R}^n,
\]
when $ f\in L^2(\mathbb{R}^n)^N$ and $F$ is a Carath\'eodory map. \eqref{1} has not been considered before.  The case of bounded domains has been studied by several authors, firstly by Campanato and under Campanato's ellipticity condition on $F$. By introducing a new much weaker notion of ellipticity, we prove solvability of \eqref{1} in a tailored Sobolev ``energy" space and a uniqueness estimate. The proof is based on the solvability of the linearised problem by Fourier transform methods, together with a ``perturbation device" which allows to use Campanato's near operators. We also discuss our hypothesis via counterexamples and give a stability theorem of strong global solutions for systems of the form \eqref{1}.
\end{abstract}

\maketitle

\section{Introduction} \label{section1}

Let $n,N\geq 2$ and let also
\[
F\ :\ \R^n \by \big(\R^N \!\ot \mS(n) \big) \larrow  \R^N
\]
be a Carath\'eodory map, namely
\[
\left\{
\begin{array}{l} 
x\mapsto F(x,\X) \text{ is measurable, for every } \X \in \R^N \!\ot \mS(n),\ms\\
\X\mapsto F(x,\X) \text{ is continuous, for almost every } x\in  \R^{n}.
\end{array}
\right.
\]
In this paper we consider the problem of existence and uniqueness of global twice weakly differentiable strong a.e.\ solutions $u :   \R^n \larrow \R^N$ to the following fully nonlinear PDE system
\beq  \label{1.1}
F(\cdot,D^2u ) \,=\,  f, \ \ \text{ a.e.\ on }\R^n,
\eeq
when $ f\in L^2(\R^n)^N$. In the above, $\mS(n)$ denotes the symmetric matrices  of $\R^{n \by n}$, $Du(x) \in \R^{N\by n}$ is the gradient matrix and $D^2u(x) \in \R^N \!\ot \mS(n) $ is the hessian tensor of $u$ at $x$. In the sequel we will employ the summation convention in repeated indices when $i,j,k,...$ run in $\{1,...,n\}$ and $\al,\be,\ga,...$ run in $\{1,...,N\}$. The standard bases of $\R^n$, $\R^N$, $\R^{N\by n}$ and $\R^N \! \ot \R^{n\by n}$ will be denoted by $\{e^i\}$, $\{e^\al\}$, $\{e^{\al i}\}$ and $\{e^{\al ij}\}$ respectively, ``$\ot$" denotes the tensor product and we abbreviate $e^{\al i} \equiv e^\al \ot e^i$,  $e^{ij} \equiv e^i \ot e^j$, $e^{\al ij} \equiv e^\al \ot e^i \ot e^j$ and $D_i\equiv \p/\p x_i$. Hence, we will write
\[
\text{$x=x_ie^i$, \ $u=u_\al e^\al$, \ $Du=(D_iu_\al) e^{\al i}$, \ $D^2u=(D^2_{ij}u_\al) e^{\al ij}$.}
\]
To the best of our knowledge, the problem \eqref{1.1} has not been considered before. However, the Dirichlet problem on bounded domains for the system $F(\cdot,D^2u )=f$ has been considered before by several authors and with different degrees of generality. The first one to address it was Campanato \cite{C1}-\cite{C4} for bounded convex $\Om$ and under a strong ellipticity condition which we recall later. Subsequent contributions to this problem and problems relevant to Campanato's work on this problem can be found in Tarsia \cite{Ta1}-\cite{Ta5}, Fattorusso-Tarsia \cite{FT1}-\cite{FT4}, Buica-Domokos \cite{BD}, Domokos \cite{Do}, Palagachev \cite{Pa1,Pa2}, Palagachev-Recke-Softova \cite{PRS}, Softova \cite{S} and Leonardi \cite{Le}. However, all vectorial contributions, even the most recent ones \cite{FT1,FT2} (wherein they consider PDE systems of the form $F(\cdot,u,Du,D^2u)=f$) are based on Campanato's original restrictive ellipticity notion, or a minor extension of it due to Tarsia \cite{Ta5}. Moreover, in the very recent papers of the author \cite{K8, K9} we are considering the relevant cases  of 1st order fully nonlinear elliptic systems $F(\cdot,Du)=f$ and also the 2nd order case $F(\cdot,D^2u)=f$ but on bounded domains.

The main consequence of Campanato's ellipticity is that the nonlinear operator $F[u]:=F(\cdot,D^2u)$ is ``near" the Laplacian $\De u$. Nearness is a functional analytic notion also introduced by Campanato in order to solve the problem, which roughly says that operators near those with ``good properties" like bijectivity inherit these properties. In the case at hand, nearness implies unique solvability of \eqref{1.1} in $(H^{2}\cap H^{1}_0)(\Om)^N$, by the unique solvability of the Poisson equation $\De u =f$ in $(H^{2}\cap H^{1}_0)(\Om)^N$ and a fixed point argument. Campanato's ellipticity relates to the Cordes condition (see Cordes \cite{Co1,Co2} and also Landis \cite{L}). 

Although Campanato's condition is stringent, it should be emphasised that in general it is not possible to obtain solvability in the class of strong solutions with the mere assumption of uniform ellipticity. Well-known counterexamples which are valid even in the linear scalar case of the second order elliptic equation 
\[
A_{ij}(x)D^2_{ij}u(x)\, =\, f(x) 
\]
with $A_{ij} \in L^\infty(\Om)$ imply that the standard uniform ellipticity $ A\geq \nu I$ does not suffice to guarantee well posedness of the Dirichlet problem when $n>2$ and more restrictive conditions are required (see e.g.\ Ladyzhenskaya-Uraltseva \cite{LU}).

In this paper we introduce a new much weaker ellipticity notion for $F$ than the Campanato-Tarsia condition and for the first time we consider the case of global solutions on $\Om=\R^n$. We prove unique solvability of \eqref{1.1} by a twice weakly differentiable map $u$ in the appropriate Sobolev space, together with a strong a priori estimate. Moreover, in the course of the proof we give a \emph{vectorial non-monotone} extension of the Miranda-Talenti inequality on the whole space. A proof of the classical Miranda-Talenti inequality in $H^2\cap H^1_0$ (Miranda \cite{M}, Talenti \cite{T}) can be found in Maugeri-Palagachev-Softova \cite{MPS}.

Our starting point for the system $F(\cdot,D^2u)=f$ is based on the analysis of the simpler case of $F$ linear in $\X$ and independent of $x$, that is when
\beq \label{1.2}
F_\al(x,\X)\, =\, \A_{\al \be i j} \X_{\be i j}.
\eeq
Here $\A$ is a linear symmetric operator $\A : \R^{N\by n}\larrow \R^{N\by n}$:
\[
\A \in \mS(N\! \by \! n), \ \  \text{ i.e. }\ \A_{\al \be ij}\, =\, \A_{\be \al j i}.
\]
For $F$ as in \eqref{1.2}, the system $F(\cdot,D^2u)=f$ becomes
\[
\A_{\al \be i j} D^2_{ij}u_\be \,=\, f_\al.
\]
By introducing the contraction operation $\A :\Z := (\A_{\al \be ij}\Z_{\al ij})e^\al$ (which extends the trace inner product $\Z:\Z=\Z_{\al i j}\Z_{\al i j}$ of $\R^N \!\ot \mS(n)$), we will write it compactly as
\beq \label{1.3}
\A: D^2u \,=\, f.
\eeq
The appropriate notion of ellipticity in this case is that the quadratic form arising from the operator $\A$
\beq \label{1.4}
\begin{array}{c}
\A \ :\ \ \R^{N\by n} \by \R^{N\by n}\ \larrow \R,\ms\\ 
\A: P\ot Q \, := \, \A_{\al \be ij}P_{\al i}Q_{\be j},
\end{array}
\eeq
is (strictly) rank-one convex on $\R^{N\by n}$, that is
\beq \label{1.5}
\A : \eta \ot a \ot \eta \ot a \, \geq \, \nu |\eta|^2|a|^2,
\eeq
for some $\nu >0$ and all $\eta \in \R^N, \ a\in \R^n$. For brevity, we will say \emph{``$\A$ is rank-one positive"} as a shorthand of the statement \emph{``the symmetric quadratic form defined by $\A$ on $\R^{N\by n}$ is rank-one strictly convex"}. Our ellipticity assumption for general $F$ is given in the following definition. We state it for a general domain $\Om\sub \R^n$:

\begin{definition}[Ellipticity] \label{def1} Let  $\Om\sub \R^n$ be open and let $F : \Om \by \big(\R^N \!\ot \mS(n) \big) \larrow  \R^N$ be a Carath\'eodory map. We call $F$ (or the system $F(\cdot,D^2u)=f$) elliptic when there exist 
\[
\begin{array}{c}
\A\,  \in \, \mS(N\!\by\! n), \text{ rank-one positive},\\
\la>\ka>0,\\
\text{$\al \in L^\infty(\Om)$, $\al>0$ a.e.\ on $\Om$ and $1/\al \in L^\infty(\Om)$,}
\end{array}
\]
such that
\beq \label{1.6}
(\A :\Z)^\top\Big[ F(x,\X+\Z) - F(x,\X)\Big] \ \geq\ \frac{\la}{\al(x)}|\A:\Z|^2\, -\ \frac{\ka}{\al(x)}\nu(\A)^2|\Z|^2,
\eeq
for all $\X,\Z \in \R^\N\ot\mS(n)$ and a.e.\ $x\in \Om \sub \R^n$. 
\end{definition}
In the above definition $\nu(\A)$ is the ellipticity constant of $\A$:
\beq \label{1.7}
\nu(\A)\, :=\, \min_{|\eta|=|a|=1} \big\{\A : \eta \ot a \ot \eta \ot a \big\}. 
\eeq
By taking 
as $\A$ the monotone tensor
\[
\A_{\al \be i j}\ =\ \de_{\al \be}\de_{ij}, 
\]
we reduce to a condition equivalent to Tarsia's notion, and by further taking $\al(x)$ constant we reduce to Campanato's notion:
\beq \label{1.8}
(\Z:I)^\top\Big[ F(x,\X+\Z) - F(x,\X)\Big] \ \geq\ c_2 |\Z:I|^2\, -\ c_1 |\Z|^2,
\eeq
$c_2>c_1>0$. In this paper, all the norms $|\cdot|$ will be the euclidean, e.g.\ on $\R^\N\ot\mS(n)$ we use $|\X|^2=\X :\X=\X_{\al ij}\X_{\al ij}$ etc, and in \eqref{1.8} we have used the obvious contraction operation $\X:X:=(\X_{\al i j}X_{ij})e^\al$. Our new ellipticity notion \eqref{1.6} relaxes \eqref{1.8} substantially: a large class of nonlinear operators to which our results apply are of the form
\[
F(x,\X)\, :=\, g^2(x)\A:\X \ + \ G(x,\X)
\]
where $\A$ rank-one positive, $g,1/g \in L^\infty(\Om)$ and $G$ is \emph{any nonlinear} map, measurable with respect to the first argument and Lipschitz with respect to the second argument, with Lipschitz constant of $G(x,\cdot)/g^2(x)$ smaller than $\nu(\A)$ (see Example \ref{ex1}). In particular, any $F \in C^{1}\big(\R^N\!\ot \mS(n)\big)^N$ such that $F'(0)$ is rank-one positive and the Lipschitz constant of $\X\mapsto F(\X)-F'(0):\X$ is smaller than $\nu(F'(0))$, is elliptic in the sense of Definition \ref{def1}. On the other hand, even if $F$ is linear, $F(\X)=\A:\X$ and in addition $\A$ defines a \emph{strictly convex} quadratic form on $\R^{N\by n}$, that is when
\[
\A : Q \ot Q \, \geq c^2 |Q|^2,\ \ \ Q\in \R^{N\by n},
\]
then $F$ may \emph{not be elliptic} in the Campanato-Tarsia sense (see Example \ref{ex2}).

The general program we deploy herein is the following: we first establish existence and uniqueness to the system \eqref{1.1} in the linear case with constant coefficients for $F(\X)=\A:\X$. Then, we use the new ellipticity notion, a ``perturbation device" which is a consequence of this ellipticity and employ Campanato's theorem of bijectivity of near operators, in order to prove existence and uniqueness for \eqref{1.1} in the general case. More precisely, in Section \ref{section3} we prove existence and uniqueness of global strong a.e.\ solutions to \eqref{1.1} in the linear case of $F(\X)=\A:\X$ when $\A$ satisifes \eqref{1.5} and $n\geq 5$. The appropriate Sobolev space is
\beq \label{1.10}
W^{2,2}_{\ast}(\R^n)^N\, :=\,   \Big\{u\in L^{2^{**}}(\R^n)^N \ \Big| \ Du \in L^{2^*}(\R^n)^{Nn},\ D^2u \in L^2(\R^n)^{Nn^2} \Big\}.
\eeq
Here the exponent $2^*$ is the conjugate exponent of $2$ and $2^{**}=(2^*)^*$:
\beq \label{1.11}
2^*\, =\, \frac{2n}{n-2}\ ,\ \ \ 2^{**}\, =\, \frac{2n}{n-4}.
\eeq
The reason why we have to restrict ourselves to dimensions $n\geq5$ relate to the Gagliardo-Nirenberg-Sobolev inequality: for $n\leq 4$, $W^{2,2}_{\ast}(\R^n)^N$ is not a Banach space with respect to the $L^2$ seminorm of the hessian. When $n\geq 5$, we prove existence, uniqueness and also an \emph{explicit} representation formula for the solution which lives in $W^{2,2}_{\ast}(\R^n)^N$ by utilising the Fourier transform (Theorem \ref{th1}). Next, in Section \ref{section4} we tackle the general case of fully nonlinear $F$ satisfying Definition \ref{def1}  (Theorem \ref{th2}). This is based on the solvability of the linear problem, our ellipticity assumption and Campanato's result of ``near operators" taken from \cite{C5}, which we recall herein for the convenience of the reader (Theorem \ref{th3}). A byproduct of our method is a strong uniqueness estimate in the form of a \emph{comparison  principle} for the distance of any solutions in terms of the distance of the right hand sides of the equations. A crucial ingredient of our analysis is the following sharp hessian estimate
\beq \label{1.12}
\big\|D^2u \big\|_{L^2(\R^n)} \, \leq \, \frac{1}{\nu(\A)}\big\|\A:D^2u \big\|_{L^2(\R^n)} 
\eeq
valid for all $u\in W^{2,2}_{\ast}(\R^n)^N$, which is established in Proposition \ref{pr2}. The inequality \eqref{1.12} is a vectorial non-monotone extension of the Miranda-Talenti inequality to the whole space and beyond the case $\A\!:\!D^2u=\De u$ of the classical result. In Section \ref{section2} we discuss some examples and counterexamples, as well as an equivalent formulation of our ellipticity condition which is the analogue of Campanato's ``A-condition". Finally, in Section \ref{section5} we discuss an extension of our main result to result of stability type for strong global solutions of fully nonlinear systems.

We note that Campanato's notion of nearness has been relaxed by Buica-Domokos in \cite{BD} to a ``weak nearness", which still retains most of the features of (strong) nearness. In the same paper, the authors also use an idea similar to ours, namely a fully nonlinear operator being ``near" a general linear operator, but they implement this idea only in the scalar case. 

Except for its intrinsic analytical interest, the motivation to study the present problem comes in part from Differential Geometry, in particular semi-Riemannian Geometry/General Relativity,but also from Conformal Geometry. In addition, understanding this problem is an important stepping stone in order to understand the problem of non-existence of minimisers in Calculus of Variations for 2nd order non-convex problems. More importantly, the present problem can be seen as a simplied version of the complicated equations which arise in the recently initiated vectorial Calculus of Variations in the space $L^\infty$ for supremal functionals. We have collected some details about these problems and how they relate to the current paper in Section \ref{Motivations_and_Applications}.

We conclude this introduction by noting that the fully nonlinear case of \eqref{1.1} has been studied also when $F$ is \emph{coercive instead of elliptic}. By using the analytic Baire category method of the Dacorogna-Marcellini \cite{DM} which is the ``geometric counterpart" of Gromov's Convex Integration, one can prove that, under certain structural and compatibility assumptions, the Dirichlet problem has \emph{infinitely many} strong a.e.\ solutions in the space $W^{2,\infty}(\R^n)^n$. However, ellipticity and coercivity of $F$ are, roughly speaking, mutually exclusive and in order to get uniqueness under this method, appropriate extra selection criteria of ``good" solutions are required, yet to be determined. 

On the other hand, the scalar theory of single elliptic equations has a much richer theory, for both classical/strong a.e.\ solutions of strongly elliptic equations, (see Gilbarg-Trudinger \cite{GT}) as well as for ``nonvariational weak solutions" of degenerate elliptic equations, namely viscosity solutions (Crandall-Ishii-Lions \cite{CIL}, Cabr\'e-Caffarelli \cite{CC} and for a pedagogical introduction see the author's monograph \cite{K}). However, except for the (fairly) broad theory for divergence strictly elliptic systems (see e.g.\ Giaquinta-Martinazzi \cite{GM}), for fully nonlinear systems the existing theory is very limited (but see the very recent development of the theory of $\mathcal{D}$-solutions in \cite{K10, K11} and Section \ref{Motivations_and_Applications}).

\section{Ellipticity, examples and counterexamples} \label{section2}

We begin by noting the simple algebraic fact that our ellipticity notion of Definition \ref{def1} implies a sort of  generalised ``non-monotone" Legendre-Hadamard condition (or strict rank-one convexity in the linear case) relative to $\A$. If $\A$ is \emph{monotone}, that is if
\[
\A_{\al \be ij}\ = \ \de_{\al \be}A_{i j},
\]
for some $A\in \mS(n)$, then we reduce to rank-one convexity. Accordingly, we have the next result:

\bl[Non-monotone rank-one convexity] \label{le1} Let $\Om \sub \R^n$ be open and let $F :  \Om \by \big(\R^{N}\!\ot \mS(n)\big) \larrow  \R^N$ be a Carath\'eodory map satisfying Definition \ref{def1} for some $\A$, $\ka,\la$ and $\al$. Then, we have the estimate
\[
\Big(F\big(x,\X+\eta \ot a \ot a\big)-F(x,\X)\Big)^\top\big(\A:\eta \ot a\ot a\big) \ \geq \, \frac{(\la-\ka)\nu(\A)^2}{\|\al\|_{L^\infty(\Om)}} |\eta|^2|a|^4,
\] 
for all $\eta \in \R^N$, $a\in \R^n$ and a.e.\ $x\in \Om$. In particular, if $F$ is linear and $F(x,\X)=\textbf{G}(x):\X$, then
\[
\big(\textbf{G}(x):\eta \ot a\ot a \big)^\top(\A:\eta \ot a\ot a) \ \geq \, \frac{(\la-\ka)\nu(\A)^2}{\|\al\|_{L^\infty(\Om)}} |\eta|^2|a|^4,
\] 
for all $\eta \in \R^N$, $a\in \R^n$. 
\el

\BPL \ref{le1}. Choose $\Z:=\eta \ot a\ot a$ for $\eta \neq 0$ and observe that
\beq \label{2.1}
|\Z|^2 \ = \ |\eta \ot a\ot a |^2\ =\eta_\al a_i a_j\, \eta_\al a_i a_j\ = \ |\eta|^2|a|^4
\eeq
and also, by \eqref{1.7}, we have
\begin{align} \label{2.2}
|\A:\Z|^2\ &= \ \max_{|\xi|=1}\, \left|\xi^\top \big(\A:\eta \ot a\ot a \big)\right|^2  \ \geq \ \left|\frac{\eta}{|\eta|}^\top \big(\A:\eta \ot a\ot a \big)\right|^2 \ = \nonumber\\
&=\ \frac{1}{|\eta|^2}\big| \eta_\al \A_{\al \be i j }\eta_\be a_i a_j \big|^2  \ \geq \frac{1}{|\eta|^2}\nu(\A)^2|\eta|^4|a|^4. 
\end{align}
Hence, by \eqref{1.6} and \eqref{2.1}, \eqref{2.2} we obtain
\begin{align}
\Big(F\big(x,\X+ \eta \ot a \ot a\big)-F(x,\X)\Big)^\top\big(\A:\eta \ot a\ot a\big) \nonumber
\end{align}
\begin{align}
 & \geq  \frac{\la}{\al(x)} \nu(\A)^2|\eta|^2|a|^4 \ - \ \frac{\ka}{\al(x)}
\nu(\A)^2|\eta|^2|a|^4 \nonumber\\
&\geq \, \frac{(\la-\ka)\nu(\A)^2}{\|\al\|_{L^\infty(\Om)}} |\eta|^2|a|^4,\nonumber
\end{align}
and the lemma ensues.       \qed

\ms

We now rewrite our ellipticity condition of Definition \ref{def1} to a formulation which is along the lines of Campanato's ``$A$-Condition" and Tarsia's ``$A_x$-Condition" (see \cite{Ta3, Ta4}).

\begin{definition}[K-Condition] \label{def2} Let $\Om \sub \R^n$ be open and $F :  \Om \by \big(\R^{N}\!\ot \mS(n)\big) \larrow  \R^N$ a Carath\'eodory map. We say  that $F$ is elliptic (or the PDE system $F(\cdot,D^2u)=f$ is elliptic) when there exist
\[
\begin{array}{c}
\A \in \mS(N\!\by\! n)\ \text{ rank-one positive},\\
\al \in L^\infty(\Om), \al>0 \text{ a.e.\ on }\Om, 1/\al \in L^\infty(\Om),\\
\be,\ga>0\ \text{ with }\be+\ga<1,
\end{array} 
\]
such that
\beq \label{2.3}
\left| \A:\Z\, -\, \al(x)\Big(F(x,\X+\Z) -F(x,\X)\Big) \right|^2 \ \leq\,  \be \nu(\A)^2|\Z|^2\ + \ \ga|\A:\Z|^2.
\eeq
\end{definition}

\ms

\noi We recall that $\nu(\A)$ is the ellipticity constant of $\A$ and is given by \eqref{1.7}. The following result certifies that the ellipticity condition of Definition \ref{def1} is equivalent to the $K$-condition of Definition \ref{def2}, if $F$ is globally Lipschitz continuous with respect to the second argument.

\begin{lemma}[Ellipticity vs $K$-Condition] \label{pr1} Let $\Om \sub \R^n$ be open and let $F :  \Om \by \big(\R^{N}\!\ot \mS(n)\big) \larrow  \R^N$ be a Carath\'eodory map. 
Then the following statements are equivalent:
\begin{enumerate}

\item There exist $\A \in \mS(N\!\by\! n)$ rank-one positive, $\be,\ga>0$ with $\be+\ga<1$ and $\al \in L^\infty(\Om)$ with $\al>0$ a.e.\ on $\Om$ and $1/\al \in L^\infty(\Om)$ with respect to which $F$ satisfies Definition \ref{def2}.

\ms

\item There exist $\A \in \mS(N\!\by\! n)$ rank-one positive, $\la>\ka>0$  and $\al \in L^\infty(\Om)$ with $\al>0$ a.e.\ on $\Om$ and $1/\al \in L^\infty(\Om)$ with respect to which $F$ satisfies Definition \ref{def1}. Moreover, $\X\mapsto F(x,\X)$ is globally Lipschitz continuous on $\R^{N}\!\ot \mS(n)$, essentially uniformly in $x\in \Om$:
\beq \label{2.4}
\underset{x\in \Om}{\ess\,\sup}\sup_{\X\neq\Y\text{ in }\R^N\!\ot \mS(n)} \frac{|F(x,\Y)-F(x,\X)|}{|\Y-\X|}\ =: \, M \ < \ \infty.
\eeq
         \end{enumerate}
\end{lemma}

\BPL \ref{pr1}. Assume (1) holds. Then, \eqref{2.3} implies
\begin{align}
|\A:\Z|^2 \ +&\ \al(x)^2 \Big|F(x,\X+\Z) -F(x,\X)\Big|^2\nonumber \\ 
&\ \ \ \ \ -\  2\al(x)(\A:\Z)^\top \Big(F(x,\X+\Z) -F(x,\X)\Big) \nonumber \\
&  \leq\ \be\nu(\A)^2|\Z|^2\,  +\, \ga|\A:\Z|^2  \nonumber
\end{align}
Hence,
\begin{align}
 \al(x)(\A:\Z)^\top &\Big(F(x,\X+\Z) -F(x,\X)\Big)  \nonumber\, \geq\,  \frac{1-\ga}{2} |\A:\Z|^2\ - \  \frac{\be}{2} \nu(\A)^2|\Z|^2 . \nonumber
\end{align}
and we obtain \eqref{1.6} for
\[
\la\ :=\ \frac{1-\ga}{2}\ , \ \ \ \ka\ :=\ \frac{\be}{2}
\]
since $\la>\ka>0$, because $\ka>0$ and $\la-\ka =\frac{1}{2}\big(1-(\be+\ga) \big)>0$. In addition, again by \eqref{2.3}, we have
\[
\al(x)\big|F(x,\X+\Z) -F(x,\X)\big| \ \leq\,  |\A:\Z|\ +\ \sqrt{\be} \nu(\A)|\Z|\ + \ \sqrt{\ga}|\A:\Z|
\]
and hence
\begin{align}
\big|F(x,\X+\Z) -F(x,\X)\big| \ &\leq\,  \frac{1}{\al(x)}\Big((1+\sqrt{\ga})|\A:\Z|\ +\ \sqrt{\be} \nu(\A)|\Z|\Big) \nonumber\\
&\leq\,  \Big\|\frac{1}{\al}\Big\|_{L^\infty(\Om)}\Big((1+\sqrt{\ga})|\A|\, +\, \sqrt{\be} \nu(\A)\Big) |\Z|.\nonumber
\end{align}
Consequently, \eqref{2.4} follows and we have just shown that (1) implies (2).
\ms

Conversely, assume (2) and fix $\si>0$. Let also $M$ be as in \eqref{2.4}. Then, by \eqref{2.4} and \eqref{1.6} we have
\[
\frac{\al(x)^2}{(\la \si)^2}\big|F(x,\X+\Z)-F(x,\X) \big|^2 \ \leq\ \left(\frac{M\al(x)}{\la \si\, \nu(\A)}\right)^2\nu(\A)^2|\Z|^2
\]
and
\begin{align}  
|\A:\Z|^2\ -\ 2\frac{\al(x)}{\la\si}(\A:\Z)&^\top \Big( F(x,\X+\Z)-F(x,\X) \Big) \nonumber\\
& \leq\ \frac{2\ka}{\la\si}\nu(\A)^2|\Z|^2\ + \Big(1-\frac{2}{\si} \Big)|\A:\Z|^2. \nonumber
\end{align}
By adding the above two inequalities, we get
\begin{align}  
&\left|\A:\Z - \frac{\al(x)}{\la\si}\Big( F(x,\X+\Z)-F(x,\X) \Big) \right|^2 \nonumber\\
& \hspace{30pt} \leq\  \left(   \frac{2\ka}{\la\si} +  \frac{1}{\si^2}\left(\frac{M\al(x)}{\la\, \nu(\A)}\right)^2 \right) \nu(\A)^2|\Z|^2 +\, \Big(1-\frac{2}{\si} \Big)|\A:\Z|^2 . \nonumber
\end{align}
Since $\ka/\la<1$, by choosing $\si>0$ large, we can arrange things such that Definition \ref{def2} is satisfied for the same $\A$ as in Definition \ref{def1} and
\[
\be\ := \   \frac{2}{\si} \left( \frac{\ka}{\la}+  \frac{1}{2\si}\left(\frac{M\|\al\|_{L^\infty(\Om)}}{\la \nu(\A)}\right)^2\right)  \ ,\ \ \ \ga\ :=\  1-\frac{2}{\si}\ , \ \ \al'(x)\, :=\ \frac{\al(x)}{\la \si},
\]
because $\be+\ga<1$, for $\si$ large. The lemma has been established.       \qed

\ms

The previous result allows us to exhibit a large class of nonlinear operators to which our existence-uniqueness results apply.

\begin{example}[A class of elliptic ``coefficients" satisfying the $K$-Condition] \label{ex1}
Nontrivial fully nonlinear examples of maps $F$ which are elliptic in the sense of the Definition \ref{def1} above are easy to find. Let $\Om\sub \R^n$ be open, $g$ measurable with $g^2, 1/g^2 \in L^\infty(\Om)$ and consider any fixed tensor $\A \in \mS(N\!\by\!n)$ for which $\nu(\A)>0$ and any Carath\'eodory map
\[
G\ : \ \Om \by \big(\R^{N}\!\ot \mS(n)\big) \larrow \R^N
\]
which is Lipschitz with respect to the second variable and
\[
\underset{x\in \Om}{\ess\,\sup}\, \left\|\frac{G(x,\cdot)}{g^2(x)}\right\|_{{Lip} \big( \R^\N\ot\mS(n) \big) } < \ \nu(\A).
\]

Then, the map $F :  \Om \by \big(\R^{N}\!\ot \mS(n)\big) \larrow \R^N$ given by
\[
F(x, \X)\, :=\, g^2(x)\A:\X\, +\, G(x,\X)
\]
satisfies Definition \ref{def2}, since there is $\be \in (0,1)$ such that
\begin{align}
\left|\A:\Z - \frac{1}{g^2(x)}\Big(F(x,\X+\Z)-F(x,\X)\Big) \right|^2\, &=\, \left|\frac{G(x,\X+\Z)-G(x,\X)}{g^2(x)} \right|^2 \nonumber\\
&\leq\, \be \,\nu(\A)^2 |\Z|^2\nonumber\\
&\leq\, \be \,\nu(\A)^2 |\Z|^2\ +\ \ga |\A:\Z|^2, \nonumber
\end{align}
and hence $F$ satisfies \eqref{2.3} for $\al:=g^{-2}$, some $\be \in (0,1)$ and any $\ga \in (\be, 1)$, e.g.\ $\ga:=(1-\be)/2$.

\noi \textit{Thus, every Lipschitz perturbation of an elliptic constant tensor gives a fully nonlinear elliptic map, when the Lipschitz constant of the perturbation is strictly smaller than the ellipticity constant of the tensor.}
\end{example}

We now show that our ellipticity condition, either in the guises of Definition \ref{def1} or in the guises of Definition \ref{def2} is \emph{strictly weaker} than the Campanato-Tarsia definition. More precisely, we give an example of a symmetric $\A \in \mS(N\by n)$ which is (not merely rank-one positive, but) positive and the respective map $F(\X):=\A:\X$ does not satisfy \eqref{1.8}. On the other hand, every such $F$ is automatically elliptic in our sense. The idea of this example is inspired by the examples in \cite{Ta4}.

\begin{example}[A strictly convex $\A$ not satisfying Campanato's $A$-Condition] \label{ex2} There exists  $\A \in \mS(2\!\by\!2)$ such that
\beq \label{2.5}
\A:Q\ot Q \ \geq\ |Q|^2, \ \  \ Q\in \R^{2\by 2},
\eeq
which is such that there do not exist constants  $c_2>c_1>0$ for which $F(\X):=\A:\X$ satisfies \eqref{1.8}.  Indeed, let us define 
\[
\A \, :=\, 
\left[
\begin{array}{c|c}
\A_{11} & \textbf{0}\\ \hline 
\textbf{0} & \A_{22}\\ 
\end{array}
\right]
\]
by setting
\[
\A_{11} \, :=\, 
I\, =\, \left[
\begin{array}{cc}
1 & 0\\
0 & 1\\ 
\end{array}
\right]\ , \ \ \ 
 \A_{22} \, :=\, 
m \left[
\begin{array}{cc}
2 & 1\\
1& 2\\ 
\end{array}
\right]\ , \ m\geq 8.
\]
In index form, this means $\A_{11ij}=\de_{ij}$, $\A_{12ij}=\A_{21ij}=0$, $\A_{2211}=\A_{2222}=2m$, $\A_{2212}=\A_{2221}=m$. Then, $\A$ satisfies \eqref{2.5}, since
\begin{align}
\A:Q\ot Q\ &=\ \A_{\al \be i j}Q_{\al i}Q_{\be j}  \nonumber\\
&=\, \A_{11ij}Q_{1i}Q_{1j} \ + \ \A_{2211}Q_{21}Q_{21}\ +\ \A_{2222}Q_{22}Q_{22}  \nonumber\\
&\ \ \ \ +\ \A_{2212}Q_{21}Q_{22}\ +\ \A_{2221}Q_{22}Q_{21}  \nonumber\\
&=\ (Q_{11})^2\ +\  (Q_{12})^2\ +\ 2m\big( (Q_{21})^2 \, +\,  (Q_{22})^2 \big)\ +\ 2mQ_{21}Q_{22}  \nonumber\\
&\geq |Q|^2\ +\ m(Q_{21}\, +\, Q_{22})^2. \nonumber
\end{align}
Suppose now that there exist $c_2>c_1>0$ such that $\X\mapsto \A:\X$ satisfies \eqref{1.8}, that is for all $\X \in \R^2\!\ot \mS(2)$,
\beq \label{2.6}
\A_{\al \be i j }\X_{\be ij}\X_{\al kk}\ \geq\ c_2\big(\X_{\be ii}\X_{\be jj}\big) \ - \ c_1\big(\X_{\al ij}\X_{\al ij} \big).
\eeq
We will show that specific choices of $\X$ lead to a contradiction and such an estimate can not hold. We first choose
\[
\X\ :=\ 
\left[
\begin{array}{c}
I \\ \hline
\textbf{0}
\end{array}
\right]
\]
that is, we take $\X_{1ij}=\de_{ij}$, $\X_{2ij}=0$. We calculate:
\begin{align}
\A_{\al \be i j }\X_{\be ij}\X_{\al kk} \ &=\ \A_{11ij}\X_{1ij}\X_{1kk}\ =\ \de_{ij}\de_{ij}\de_{kk}\ = \ 4, \nonumber\\
\X_{\be ii}\X_{\be jj}\ &= \ \X_{1 ii}\X_{1 jj}\ = \ \de_{ii}\de_{jj}\ = \ 4,  \nonumber\\
\X_{\al ij}\X_{\al ij} \ &=\ \X_{1 ij}\X_{1ij} \ =\ \de_{ij} \de_{ij}\ =\ 2. \nonumber
\end{align}
Then, \eqref{2.6} implies $4\geq 4c_2 - 2c_1$, and since $c_1<c_2$, we obtain
\beq \label{2.7}
c_2\ <\ 2.
\eeq
Next, we choose
\[
\X\ :=\ 
\left[
\begin{array}{cc}
 & \! \! \! \! \! \! \! \! \! \! \! \!\textbf{0} \\ \hline
-1 & \ 3\\
\ 3& -1
\end{array}
\right]
\]
that is, we take $\X_{1ij}=0$, $\X_{211}=\X_{222}=-1$, $\X_{212}=\X_{221}=3$. We calculate:
\begin{align}
\A_{\al \be i j }\X_{\be ij}\X_{\al kk} \ &=\ \A_{22ij}\X_{2ij}\X_{2kk}\ =\ -2m(-2-2+3+3)\ = \ -4m, \nonumber\\
\X_{\be ii}\X_{\be jj}\ &= \ \X_{2 ii}\X_{2 jj}\ = \ (\X_{211}\, +\, \X_{222})^2\ = \ 4,  \nonumber\\
\X_{\al ij}\X_{\al ij} \ &=\ \X_{2 ij}\X_{2ij} \ =\ (-1)^2+(-1)^2+3^2+3^2\ =\ 20. \nonumber
\end{align}
Hence, \eqref{2.6} implies $-4m\geq4c_2-20c_1$, which by using that $c_1<c_2$ gives
\beq \label{2.8}
c_2\ >\ \frac{m}{4}.
\eeq
Since \eqref{2.7} and \eqref{2.8} are incompatible, we see that \eqref{1.8} can not be satisfied by $\A$ and $\X\mapsto \A:\X$ is not elliptic in the Campanato sense. However, $\A$ is automatically elliptic in our sense since it satisfies Definition \ref{def1}.
\end{example}

Tarsia proved in \cite{Ta5} that in the scalar case of $N=1$ and for $F(x,\X)$ linear, that it when $F(x,\X)=\A(x):\X$, the ``$A_x$-Condition" (namely \eqref{2.3} for $N=1$ and $\A_{\al \be ij}=\de_{ij}$) is equivalent to a condition with perhaps different $\be$ and $\al$, \emph{but with $\ga=0$}. In other words, in the linear case, the term of the trace $|\Z:I|$ can be absorbed into the term of the norm $|\Z|$. This result has been simplified by Domokos in \cite{D}. Now we show that in the nonlinear case this is not in general possible, not even in the scalar case. Hence, we deduce that our $K$-Condition of Definition \ref{def2} \emph{can not} be simplified to a condition with $\ga=0$.

\begin{example}[Optimality of the $K$-Condition] \label{ex3} For any $n\geq 3$, there exists a Lipschitz function $F\in C^{0,1}\big(\mS(n)\big)$ which satisfies 
\beq \label{2.9}
\Big| Z:I \, -\, \al\Big(F(X+Z)-F(X)\Big) \Big|^2 \ \leq\ \be |Z|^2\ +\ \ga |Z:I|^2
\eeq
for some $\al,\be,\ga>0$ with $\be+\ga<1$ and all $X,Z \in \mS(n)$, but \emph{does not satisfy} \eqref{2.9} with $\ga=0$ for any $\be \in (0,1)$.  The appropriate $F :\mS(n)\larrow \R$ is given by
\beq \label{2.10}
F(X)\ :=\ X:I \, -\, b|X|\,-\,c|X:I|
\eeq
where the parameters $b,c$ satisfy
\[
c \, >\, b\, >\,0,\ \ \sqrt{n}c\, +\, b\, >\, 1,\ \ b\, +\, c\, <\, 1, \ \ b^2+\,c^2<\,\frac{1}{2}.
\]
A specific choice of such values is
\[
c\, =\, \frac{1}{\sqrt{n}}\ , \ \ \ b\, =\, \frac{1}{10+\sqrt{n}}.
\]
We begin with the next claim:
\begin{claim} Let $F$ be given by \eqref{2.10}. Then, $F$ satisfies \eqref{2.9} with
\[
\ga(\al) \, :=\, 2\big(|1-\al|\, +\, \al c\big)^2\ ,\ \ \ \ \be(\al) \, :=\, 2(\al b)^2, 
\]
if and only if $\al \in (1-\al_0,1+\al_0)$ for some $\al_0 =\al_0(b,c) \in (0,1)$.
\end{claim}

\bp For any $\al>0$, we have
\begin{align}
\Big| Z:I \, -\, \al\Big(F(X+Z)-F(X)\Big) \Big|^2 \ & =\ \Big| Z:I -\al Z:I +\al 
\Big[b\big(|X+Z|-|X| \big)  \nonumber \\
&\ \ \ \ \ + c\big(|X:I+Z:I|-|X:I| \big) \Big] \Big|^2  \nonumber\\
& \leq\ \Big(\big( |1-\al|\, +\, \al c \big)|Z:I|\ +\ \al b |Z|\Big)^2  \nonumber\\
& \leq\ 2\big( |1-\al|\, +\, \al c \big)^2|Z:I|^2\ +\ 2(\al b)^2 |Z|^2  \nonumber\\
&=\ \ga(\al) |Z:I|^2 \ +\ \be(\al) |Z|^2. \nonumber
\end{align}
Now note that
\[
\ga(\al) \, +\, \be(\al)\, =\, 2\big(|1-\al|\, +\, \al c\big)^2 +\, 2(\al b)^2\ \larrow\ 2(c^2 +\,b^2)\, <\, 1,
\]
as $\al \ri 1$. Hence, there is an $\al_0>0$ such that $\ga(\al)+ \be(\al)<1$ if and only if $|\al-1|<\al_0$. Conversely, if $|\al-1|\geq\al_0$, then $\ga(\al)+ \be(\al)\geq1$. We choose $X_0:=0$ and $Z_0$ such that
\[
\sqrt{\frac{\be(\al)}{2}} |Z_0| \ =\ \sqrt{\frac{\ga(\al)}{2}} |Z_0:I|,\ \ \ (1-\al)Z_0:I\,>\,0.
\]
As explicit such $Z_0$ is
\[
Z_0\, :=\, 
\left[
\begin{array}{cc|c}
\sgn(1-\al) & t  & 0\\
t & 0 & 0\\ 
\hline
0  &  0 & \textbf{0} 
 \\
\end{array}
\right], \ \ \ t\,:=\, \sqrt{\frac{1}{2}\left(\frac{\ga(\al)}{\be(\al)}-1\right)}
\]
which is admissible choice since
\[
\frac{\ga(\al)}{\be(\al)} \, \geq\, \frac{c^2}{b^2}\, >\, 1.
\]
Then we have
\begin{align}
\Big| Z_0:I \, -\, \al\Big(F(X_0+Z_0)-F(X_0)\Big) \Big|^2 \ & =\ \Big| Z_0:I -\al Z_0:I +\al b |Z_0|  +\al c | Z_0:I|  \Big|^2  \nonumber\\
& =\ \Big| (1-\al) Z_0:I +\al b |Z_0|  +\al c | Z_0:I|  \Big|^2  \nonumber\\
& =\ \Big( \big(|1-\al|+\al c\big) |Z_0:I|  +\al b | Z_0|  \Big)^2  \nonumber\\
& =\ \left(\sqrt{\frac{\ga(\al)}{2}} |Z_0:I|\, +\, \sqrt{\frac{\be(\al)}{2}} |Z_0| \right)^2 \nonumber\\
&=\ \ga(\al) |Z_0:I|^2 \ +\ \be(\al) |Z_0|^2 \nonumber
\end{align}
and hence the estimate does not hold if $|\al-1|\geq\al_0$.
\ep

We now have the next:

\begin{claim} Let $F$ be given by \eqref{2.10}. Then, for any $\al \in (1-\al_0,1+\al_0)$, there exist $Z^0,X^0 \in \mS(n)$ with $Z^0\neq 0$ such that
\beq \label{2.11}
\Big| Z^0:I \, -\, \al\Big(F(X^0+Z^0)-F(X^0)\Big) \Big| \  =\ |Z^0|.
\eeq
Hence, the estimate \eqref{2.9} can not hold with $\ga=0$ for any $\be\in (0,1)$, regardless the choice of admissible $\al$.
\end{claim}

\bp For any fixed $\al \in (1-\al_0,1+\al_0)$, we choose 
\[
X^0\, :=\, \textbf{0}\ , \ \ \  Z^\pm\, :=\, 
\left[
\begin{array}{cccc}
\pm \zeta & 1 & \cdots & 1 \\
1 & \ddots & &  1\\
\vdots & & \ddots & \vdots\\
1 & \cdots & 1& \pm \zeta
\end{array}
\right]
\ ,\ \ \ 
Z^0\, :=\, 
\left\{
\begin{array}{l}
Z^+,\ \ \text{ when } 1-\al_0<\al\leq 1,\\
Z^-,\ \ \text{ when } 1<\al<1+\al_0,
\end{array}
\right.
\]
and
\[
\zeta\ :=\ \frac{(1-\al b) \sqrt{n-1}}{\sqrt{ n\big(|1-\al| + \al c\big)^2\, -\, \big(1-\al b\big)^2}}.
\]
Since we have chosen $b,c$ such that $\sqrt{n}c+b>1$,  $\zeta$ is well defined as a real number: indeed, by elementary algebra, we have
\[
n\big(|1-\al| + \al c\big)^2\, -\, \big(1-\al b\big)^2>\, 0 \ \ \ \Longleftrightarrow\ \ \ \sqrt{n}c\, +\, b\, >\, \sqrt{n}\, -\, \frac{\sqrt{n}-1}{\al}\, =:\, \chi(\al),
\]
when $1-\al_0<\al\leq 1$, and 
\[
n\big(|1-\al| + \al c\big)^2\, -\, \big(1-\al b\big)^2>\, 0 \ \ \ \Longleftrightarrow \ \ \ \sqrt{n}c\, +\, b\, >\,  \frac{1+\sqrt{n}}{\al}\, -\, \sqrt{n}\, =:\, \psi(\al),
\]
when $1\leq \al<1+\al_0$. Since both $\chi(\al)$ and $\psi(\al)$ are maximised when $\al=1$ and $\chi(1)=\psi(1)=1$, we deduce that indeed $\zeta \in \R$. We now show that $F$ satisfies \eqref{2.11} for these choices of $X^0,Z^0$. Indeed, we have that
\[
(1-\al) (Z^0:I)\ =\ |1-\al|\, \big|Z^0:I \big|
\]
and hence
\begin{align}
\Big| Z^0:I \, -\, \al\Big(F(X^0+Z^0)-F(X^0)\Big) \Big| \ & =\ \Big| (1-\al) Z^0:I\ +\ \al b|Z^0|\  +\ \al c|Z^0:I| \Big|  \nonumber\\
 & =\ \big(|1-\al|  \, +\,  \al c \big)|Z^0:I|\ +\ \al b|Z^0|    \nonumber\\
&=:\ \tilde{\ga}(\al) |Z^0:I| \ +\ \tilde{\be}(\al) |Z^0|. \nonumber
\end{align}
We conclude by showing that $Z^0$ (and any multiple of it) solves the algebraic equation
\beq \label{2.12}
\tilde{\ga}(\al) |Z^0:I| \ +\ \tilde{\be}(\al) |Z^0|\ =\ |Z^0|.
\eeq
By the definition of $\zeta$, we have
\[
\zeta\ =\ \frac{\big(1-\tilde{\be}(\al)\big) \sqrt{n-1}}{\sqrt{ n\tilde{\ga}^2(\al)\, -\, \big(1-\tilde{\be}(\al)\big)^2}}
\]
and by the definition of $Z^0$, we have
\[
|Z^0:I|\ =\ n\zeta\ , \ \ \ |Z^0|^2\, =\ n^2\, -\, n\, +\, n\zeta^2.
\]
Thus, we have
\begin{align}
\Pi(\al)\ &:=\ \big(1-\tilde{\be}(\al)\big)^2|Z^0|^2\ -\ \tilde{\ga}^2(\al)|Z^0|^2 \nonumber\\ 
&=\ \big(1-\tilde{\be}(\al)\big)^2\Big( n(n-1)\, +\, n\zeta^2 \Big)\ -\  \tilde{\ga}^2(\al)n^2\zeta^2 \nonumber\\
&= \ n\Bigg\{ 
 \big(1-\tilde{\be}(\al)\big)^2(n-1)\ +\ \big(1-\tilde{\be}(\al)\big)^2\frac{ \big(1-\tilde{\be}(\al)\big)^2(n-1)}{n\tilde{\ga}^2(\al)\, -\, \big(1-\tilde{\be}(\al)\big)^2} \nonumber\\
&\ \ \ \ \  \ \ \ -\  \tilde{\ga}^2(\al) n\frac{ \big(1-\tilde{\be}(\al)\big)^2(n-1)}{n\tilde{\ga}^2(\al) \, -\, \big(1-\tilde{\be}(\al)\big)^2}\Bigg\} \nonumber\\
&=\ 0. \nonumber
\end{align}
The conclusion follows by observing that the equation $\Pi(\al)=0$ is equivalent to \eqref{2.12}.
\ep
\end{example}

\section{Existence-uniqueness-representation for linear systems} \label{section3}

In this section we prove unique solvability in the case of the linear system
\beq  \label{4.1}
\A:D^2u \,=\,  f, \ \ \text{ a.e.\ on }\R^n,
\eeq
by a global solution in $W^{2,2}_\ast(\Om)^N$ for any $ f\in L^2(\R^n)^N$, when $\A\in \mS(N\by n)$ is strictly rank-one positive and $n\geq 5$. The functional ``energy" space $W^{2,2}_\ast(\Om)^N$ is given by  \eqref{1.10}. We note that in \eqref{1.10} the meaning of ``$L^{2^*}$, $L^{2^{**}}$" is $``L^p$ for $p=2^*, 2^{**}$" and not the dual or bidual space. The exponents $2^*$ and $2^{**}$ are given by \eqref{1.11}. The elementary ideas of Fourier Analysis we use herein can be found e.g.\ in Folland \cite{F} and we follow more or less the same notations as therein. In particular, for the Fourier transform and its inverse we use the conventions
\[
\widehat{u}(z)\,= \int_{\R^n}u(x)e^{-2\pi i x\cdot z}dx\ , \ \ \ \overset{\vee}{u}(x)\,= \int_{\R^n}u(z)e^{2\pi i x\cdot z}dz.
\]
Here ``$\cdot$" is the inner product of $\R^n$. Moreover, it is easy to see that if $\A\in \mS(N\!\by  n)$, then, in view of \eqref{1.7}
\beq \label{4.2}
\nu(\A)>0 \ \ \ \Longleftrightarrow\ \ \ \min_{|a|=1}\big\{ \det(\A : a\ot a) \big\} \, >\, 0,
\eeq
where $\A :a \ot a$ is the symmetric $N \!\by\! N$ matrix
\[
\A :a\ot a \, :=\, (\A_{\al \be i j}\, a_i a_j )\, e^\al \ot e^\be \ \in \ \mS(N).
\]
With ``sgn" we will denote the sign function on $\R^n$, namely $\sgn(x)=x/|x|$ when $x\neq 0$ and $\sgn(0)=0$. With ``$\cof(X)$" we will denote the cofactor matrix of $X \in \R^{N \by N}$ and we will tacitly use the identity
\[
X\cof(X)^\top \ =\ \cof(X)^\top X \ = \ \det(X)I. 
\]
The following are the two main results of this section. Proposition \ref{pr2} below is a variant of the Miranda-Talenti lemma from the case of the Laplacian (see e.g. \cite{MPS}) to the case of general $\A$ and on the whole space.

\begin{proposition}[The hessian estimate in $W^{2,2}_\ast(\R^n)^N$] \label{pr2}

Let  Let $n\geq 5$, $N\geq 2$ and $\A\in \mS(N\by n)$ rank-one positive with ellipticity constant $\nu(\A)$ given by \eqref{1.7}. Then, we have the estimate
\beq \label{4.3}
\big\|D^2u \big\|_{L^2(\R^n)} \, \leq \, \frac{1}{\nu(\A)}\big\|\A:D^2u \big\|_{L^2(\R^n)} 
\eeq
valid for all $u\in W^{2,2}_{\ast}(\R^n)^N$ (the space is given by \eqref{1.10}).
\end{proposition}

\bt[Existence-Uniqueness-Representation for the linear problem] \label{th1} Let $n\geq 5$, $N\geq 2$ and $\A \in \mS(N\!\by\! n)$ a rank-one positive tensor. Let also $f\in L^2(\R^n)^N$. Then, the problem
\[
\A:D^2u \,=\,  f, \ \ \text{ a.e.\ on }\R^n,
\]
has a unique solution $u$ in the space $W^{2,2}_\ast(\R^n)^N$ (given by \eqref{1.10}), which satisfies the estimate
\beq \label{4.4a}
\|u\|_{L^{2^{**}}(\R^n)}\, +\|Du\|_{L^{2^{*}}(\R^n)}\, +\|D^2u\|_{L^{2}(\R^n)}\ \leq \ C\|f\|_{L^2(\R^n)}
\eeq
for some $C>0$ depending only on $\emph{\A}$ and the dimensions, and also satisfies the estimate \eqref{4.3}. 

Moreover, we have the following representation formula for the solution:
\beq \label{4.4b}
u\, =\, -\frac{1}{4\pi^2} \lim_{m\ri \infty}
\left\{  \widehat{h_m} \ast
\left[ 
\frac{\  \emph{\cof} \,(\emph{\A} :\sgn \ot \sgn )^\top}{ \det (\emph{\A}: \sgn \ot \sgn) \ } \overset{\vee}{f}
\right]^{\wedge} 
\right\}.
\eeq
In \eqref{4.4b} $(h_m)^\infty_1 \sub \S(\R^n)$ is any sequence of even functions in the Schwartz class for which 
\[
\text{$0\, \leq\, h_m(x) \,\leq\, \frac{1}{|x|^2}$ \ \, and \, \ $h_m(x) \larrow \frac{1}{|x|^2}$, \ for a.e. $x\in \R^n$,\ \,  as $m\ri \infty$.} 
\]
The limit in \eqref{4.4b} is meant in the weak $L^{2^{**}}$ sense as well as a.e.\ on $\R^n$, and $u$ is independent of the choice of sequence $(h_m)^\infty_1 $.

\et

\begin{remark} The solution $u$ in \eqref{4.4b} is vectorial but \emph{real}, although the formula \eqref{4.4b} involves complex quantities. 
\end{remark}

\begin{remark}[Equivalent norms on $W^{2,2}_\ast(\R^n)^N$] \label{rem} When $n\geq 5$, the Gagliardo-Nirenberg-Sobolev inequality (see e.g.\ Evans \cite{E})
\[
\|v\|_{L^{p^*}(\R^n)} \ \leq C(n,p)\|Dv\|_{L^p(\R^n)}, \ \ \ p^*\, =\, \frac{np}{n-p},
\]
applied to $Du$ for $p=2$ and to $u$ for $p=2^*$, imply that two equivalent norm on $W^{2,2}_\ast(\R^n)^N$ are
\[
\|D^2u\|_{L^2(\R^n)}\ \approx\ \|u\|_{W^{2,2}_\ast(\R^n)}\ :=\ \|u\|_{L^{2^{**}}(\R^n)}\, +\, \|Du\|_{L^{2^*}(\R^n)}\, +\, \|D^2u\|_{L^2(\R^n)}. 
\]
\end{remark} 

The first step towards the hessian estimate is the next simple algebraic lemma, which allows to use Plancherel's theorem. 

\bl[Extension of rank-one convexity on $\C^{N\by n}$] \label{le2} Let $\A\in \mS(N\by n)$ be rank-one positive, that is
\[
\A:\eta \ot a \ot \eta \ot a\ \geq \nu |\eta|^2|a|^2, \ \ \ \eta\in \R^N,\ a\in \R^n. 
\]
We extend the quadratic form arising from $\A$ as a Hermitian form on $\C^{N\by n}$ by setting
\[
\A\ :\ \C^{N\by n} \by \C^{N\by n}\larrow \C,\ \ (P,Q)\mapsto \A: P\ot \overline{Q}.
\]
Then, we have that $\A: Q\ot \overline{Q} \in \R$ and also
\[
\A:\xi \ot a \ot \overline{\xi} \ot a\ \geq \nu |\xi|^2|a|^2, \ \ \ \xi\in \C^N,\ a\in \R^n. 
\]
\el
\noi We note that the norms on the complex spaces are the euclidean: $|\xi|^2=\xi_{\al} \overline{\xi_\al}$, etc.

\BPL \ref{le2}. The arguments are elementary, but we give them for completeness. By the symmetry of $\A$, we have
\[
\A: Q\ot \overline{Q}\ =\ \A_{\al \be i j} Q_{\al i}\overline{Q_{\be j}}\ =\ \A_{\be \al ji}\overline{Q_{\be j}} Q_{\al i}\ =\ \overline{\A_{\be \al ji}Q_{\be j}\overline{Q_{\al i}}} \ = \ \overline{\A: Q\ot \overline{Q}}.
\]
Hence, $\A: Q\ot \overline{Q} \in \R$. Next, we split $\C^N \ni \xi =\eta+i\theta$ and use symmetry again to calculate 
\begin{align}
\A:\xi \ot a \ot \overline{\xi} \ot a\ &=\  \A:(\eta+i\theta) \ot a \ot (\eta-i\theta) \ot a  \nonumber\\
&=\ \A:\eta \ot a \ot \eta \ot a\ -\ i \A:\eta \ot a \ot \theta \ot a  \nonumber\\
&\ \ \ \ +i \A:\theta \ot a \ot \eta \ot a\ +\ \A:\theta \ot a \ot \theta \ot a  \nonumber\\
&= \A:\eta \ot a \ot \eta \ot a\ +\ \A:\theta \ot a \ot \theta \ot a  \nonumber\\
&\geq \ \nu\big(|\eta|^2|a|^2 \ +\ |\theta|^2|a|^2\big)  \nonumber\\
&=\ \nu |\xi|^2|a|^2. \nonumber
\end{align}
Hence, the lemma ensues.             \qed

\ms

\BPP \ref{pr2}. We will prove the estimate when $u\in C^{\infty}_c(\R^n)^N$. In view of Remark \ref{rem}, the general case follows by a standard approximation argument. Given such a $u$, we set
\beq \label{3.2}
\A_{\al \be kl}D^2_{kl}u_\be\ =:\ f_\al \ \in \  C^\infty_c(\R^n).
\eeq
By applying the Fourier transform to the above equality, we have
\[
\A_{\al \be kl}\widehat{D^2_{kl}u_\be}\ =\ \widehat{f_\al} 
\]
and hence
\[
\A_{\al \be kl}\widehat{u_\be}(z)(2\pi iz_k) (2\pi iz_l)\ = \ \widehat{f_\al}(z) 
\]
for a.e.\ $z\in \R^n$. By multiplying by $\overline{\widehat{u_\al}}(z)$ and summing in $\al$, we get
\[
-4\pi^2\A_{\al \be kl}\widehat{u_\be}(z)z_l\overline{\widehat{u_\al}}(z) z_k\ =\ \widehat{f_\al}(z) \overline{\widehat{u_\al}}(z),
\]
a.e.\ on $\R^n$. We rewrite it as
\beq \label{3.3}
4\pi^2\A : \widehat{u}(z)\ot z\ot \overline{\widehat{u}(z)} \ot z\ =\ -\widehat{f}(z) \cdot \overline{\widehat{u}}(z).
\eeq
By Lemma \ref{le2}, both sides of \eqref{3.3} are real and positive, and also (in view of \eqref{1.7}) \eqref{3.3} implies
\beq  \label{3.4}
4\pi^2|\widehat{u}(z)|^2|z|^2\ \leq \ -\frac{1}{\nu(\A)}\widehat{f}(z) \cdot \overline{\widehat{u}}(z).
\eeq
Now we calculate:
\begin{align}
\big|\widehat{D^2u}(z) \big|^2\ &=\ \big|\widehat{u}(z)\ot (2\pi iz)\ot  (2\pi iz) \big|^2 \nonumber\\
&=\ \big|\widehat{u}(z)\ot (2\pi iz)\big|^2\, |2\pi iz|^2 \nonumber\\
&=\ 4\pi^2|\widehat{u}(z)|^2|z|^2 \, |2\pi iz|^2, \nonumber
\end{align}
for a.e.\ $z\in \R^n$. In view of \eqref{3.4}, we obtain the estimate
\begin{align}
\big|\widehat{D^2u}(z) \big|^2\ & \leq\ -\frac{1}{\nu(\A)} \widehat{f}(z) \cdot \overline{\widehat{u}}(z)\, |2\pi iz|^2 \nonumber\\
&\leq\ \frac{1}{\nu(\A)} \big|\widehat{f}(z)\big|\, \big|\overline{\widehat{u}}(z)\big| \, |2\pi iz|^2 \nonumber\\
&=\ \frac{1}{\nu(\A)} \big|\widehat{f}(z)\big|\, \big|\widehat{u}(z)\ot (2\pi iz)\ot  (2\pi iz) \big| \nonumber\\
&=\ \frac{1}{\nu(\A)} \big|\widehat{f}(z)\big|\, \big|\widehat{D^2u}(z) \big|, \nonumber
\end{align}
for a.e.\ $z\in \R^n$. For $\e\in (0,\nu(\A))$ and by Young's inequality, the above estimate gives
\[
\big|\widehat{D^2u} \big|^2\  \leq \ \frac{1}{\nu(\A)} \Big( 
\frac{1}{4\e} \big|\widehat{f}\big|^2\ +\ \e\big|\widehat{D^2u}  \big|^2    \Big),
\]
a.e.\ on $\R^n$, which, in view of \eqref{3.2}, we rewrite as
\[
\big|\widehat{D^2u} \big|^2\  \leq \ \frac{1}{4\e\big(\nu(\A) -\e\big)}  \big|\A:\widehat{D^2u}\big|^2.
\]
We choose  $\e := {\nu(\A)}/{2}$ which is the choice which maximises the denominator of the above inequality giving the value $\nu(\A)^2$, and integrate oven $\R^n$, to obtain
\[
\big\|\widehat{D^2u} \big\|^2_{L^2(\R^n)} \, \leq \, \frac{1}{\nu(\A)^2}\big\|\A:\widehat{D^2u} \big\|^2_{L^2(\R^n)}. 
\]
By applying Plancherel's theorem, the desired estimate ensues.  \qed 

\ms

\noi \textbf{Formal derivation of the representation formula.} Before giving the rigorous proof of Theorem \ref{th1}, it is very instructive to derive \emph{formally} a representation formula for the solution of $\A  : D^2u =f$. By applying the Fourier transform to the PDE system, we have
\[
\A : \widehat{D^2u}\,=\, \widehat{f}, \ \ \text{ a.e.\ on }\R^n,
\]
and hence,
\[
-4\pi^2 \, \A : \widehat{u}(z) \ot z \ot z \,= \, \widehat{f}(z),  \ \  \text{ for a.e. }z\in \R^n.
\]
For clarity, let us also rewrite this equation in index form:
\[
\left(\A_{\al \be ij} {z_i}{z_j}\right)  \widehat{u_\be}(z)\,=\, -\frac{1}{4\pi^2}\widehat{f_\al}(z).
\]
Hence, we have
\[
\left(\A :\frac{z}{|z|}\ot \frac{z}{|z|}\right)  \widehat{u}(z) \, = \, -\frac{1}{4\pi^2 |z|^2}\widehat{f}(z)
\]
and by using the identity (see \eqref{4.2})
\beq \label{iden}
\big(\A :\sgn(z)\ot \sgn(z)\big)^{-1}\, =\, \frac{\ \cof \big(\A :\sgn(z)\ot \sgn(z) \big)^\top}{ \det \big( \A :\sgn(z)\ot \sgn(z) \big)}
\eeq
we get
\begin{align}
\widehat{u}(z) \, &=\,  -\frac{1}{4\pi^2 |z|^2}\big( \A :\sgn(z)\ot \sgn(z)\big)^{-1}  \widehat{f}(z)  \nonumber\\
&=\, -\frac{1}{4\pi^2 |z|^2} \frac{\ \cof \big(\A :\sgn(z)\ot \sgn(z) \big)^\top}{ \det \big( \A :\sgn(z)\ot \sgn(z)\big)}\widehat{f}(z). \nonumber
\end{align}
By the Fourier inversion formula and  the identity $f^{\vee}(z)=\widehat{f}(-z)$, we obtain
\begin{align}
u \,&= \, -\frac{1}{4\pi^2} \left\{\frac{1}{|\cdot|^2}\frac{\ \cof \big( \A :\sgn \ot \sgn \big)^\top}{ \det \big(\A :\sgn\ot \sgn \big)} \widehat{f} \right\}^{\vee}    \nonumber\\
&= \, -\frac{1}{4\pi^2} \left\{\frac{1}{|\cdot|^2}\frac{\ \cof \big( \A :\sgn \ot \sgn \big)^\top}{ \det \big(\A :\sgn(z)\ot \sgn(z)  \big)}\overset{\vee}{f}\right\}^{\wedge} . \nonumber
\end{align}
Hence, we get the formula
\beq \label{4.5}
u \,= \,-\frac{1}{4\pi^2}  \widehat{ \frac{1}{|\cdot|^2}} \ast \left[ \frac{\ \cof \big( \A :\sgn \ot \sgn \big)^\top}{ \det \big( \A :\sgn \ot \sgn \big)}\overset{\vee}{f}\right]^{\wedge} .
\eeq
Formula \eqref{4.5} is ``the same" as \eqref{4.4b}, \emph{if we are able to pass the limit inside the integrals} of the convolution and the Fourier transform. However, in general this may not be possible. Convergence needs to be rigorously justified, and this is part of the proof of Theorem \ref{th1}. Further, by using the next identity (which follows by the properties of the Riesz potential)
\[
\left(\frac{1}{|\cdot|^2}\right)^{\wedge} =\, \ga_{n-2} \frac{1}{\ |\cdot|^{n-2}}
\]
where the constant $\ga_{\al}$ equals 
\[
\ga_\al \,=\, \frac{2^\al\, \pi^{n/2} \,\Gamma(\al/2)}{\Gamma(n/2-\al/2)}, \ \ \  \, 0\, <\, \al\, <\, n,
\]
we may rewrite \eqref{4.5} as
\beq \label{4.6}
u \, =\, -\frac{\ga_{n-2}}{4\pi^2 |\cdot|^{n-2}} \ast \left[ \frac{\ \cof \big(\A : \sgn \ot \sgn\big)^\top}{ \det \big(\A :\sgn \ot \sgn\big)}\overset{\vee}{f}\right]^{\wedge}.
\eeq
Formula \eqref{4.6} is the formal interpretation of the expression \eqref{4.4b}, which we will now establish rigorously.

\BPT \ref{th1}. By Proposition \ref{pr2}, we have the a priori estimate \eqref{4.3} for the solution, so it remains to prove existence of $u$ and the desired formula \eqref{4.4b}. Let $(h_m)^\infty_1 \sub \S(\R^n)$ be any sequence of even functions in the Schwartz class for which 
\beq \label{4.7a}
\text{$0\leq h_m(x) \leq \frac{1}{|x|^2}$ \  and \ $h_m(x) \larrow \frac{1}{|x|^2}$, \ for a.e.\ $x\in \R^n$,\  as $m\ri \infty$.} 
\eeq
We set:
\beq \label{4.8}
u_m\, :=\, -\frac{1}{4\pi^2} 
 \widehat{h_m} \ast
\left[ 
\frac{ \ \cof \,( \A : \sgn \ot \sgn )^\top }{ \det (\A : \sgn \ot \sgn) }\, \overset{\vee}{f}
\right]^{\wedge} .
\eeq
We will now show that the function $ u_m$ of \eqref{4.8} satisfies
\[
u_m\, \in\, \bigcap_{2\leq r \leq \infty}L^r(\R^n)^N\bigcap C^\infty(\R^n)^N. 
\]
Indeed, observe first that since ${h_m} \in \S(\R^n)$ and the Fourier transform is bijective on the Schwartz class, we have
\[
\widehat{h_m}\, \in \,\S(\R^n) \, \sub \, L^1(\R^n)\cap L^2(\R^n).
\]
Let now $p\in [1,2]$ and define $r$ by
\[
r\ :=\ \frac{2p}{2-p}.
\]
Then, we have
\[
1\, +\, \frac{1}{r}\, =\, \frac{1}{p}\, +\, \frac{1}{2}, \ \ \ 1\leq p\leq2,
\]
and by Young's inequality and Plancherel's theorem, we obtain
\begin{align}
\|u_m\|_{L^r(\R^n)}\ &\leq \  \frac{1}{4\pi^2} \big\| \widehat{h_m} \big\|_{L^p(\R^n)} 
\left\| \left[ 
\frac{ \ \cof \,\big( \A : \sgn \ot \sgn \big)^\top }{ \det \big( \A : \sgn \ot \sgn \big) }\, \overset{\vee}{f}
\right]^{\wedge}   
\right\|_{L^2(\R^n)}   \nonumber\\
 &\leq \  \frac{1}{4\pi^2} \big\| \widehat{h_m} \big\|_{L^p(\R^n)} 
\left\|
\frac{ \ \cof \,\big( \A : \sgn \ot \sgn \big)^\top }{ \det \big( \A : \sgn \ot \sgn \big) }\, \overset{\vee}{f}
\right\|_{L^2(\R^n)} \nonumber .
\end{align}
We now recall that the estimate \eqref{2.2}  implies 
\[
\underset{z\in \R^n}{\ess\,\inf}\, \big|\det (\A : \sgn(z) \ot \sgn(z))\big|\, >\, 0 
\]
and hence we get
\begin{align}
\|u_m\|_{L^r(\R^n)}\ 
&\leq \ \frac{1}{4\pi^2} \big\| \widehat{h_m} \big\|_{L^p(\R^n)} 
\left\|
\frac{ \, \cof \,\big( \A : \sgn \ot \sgn \big) }{ \det \big( \A : \sgn \ot \sgn \big) }\right\|_{L^\infty(\R^n)}  
\big\|\overset{\vee}{f}
\big\|_{L^2(\R^n)}   \nonumber\\
&\leq \  C \big\| \widehat{h_m} \big\|_{L^p(\R^n)} 
\left\|f
\right\|_{L^2(\R^n)}, \nonumber
\end{align}
for some $C>0$ depending only on $|\A|$ and $\nu(\A)$. Consequently, $u_m  \in L^r(\R^n)^N$ for all $r \in [2, \infty]$. Moreover, since $\widehat{h_m}\in \S(\R^N)$, we have that $u_m \in C^\infty(\R^n)^N$ by the properties of convolution. 

Next, by \eqref{4.8} and the properties of convolution, we obtain
\[
u_m\, =\, -\frac{1}{4\pi^2} 
\left[ 
h_m \frac{ \ \cof \,\big( \A : \sgn \ot \sgn \big)^\top }{ \det \big( \A : \sgn \ot \sgn \big) }\, \overset{\vee}{f}
\right]^{\wedge} ,
\]
on $\R^n$. The Fourier inversion theorem gives
\[
\overset{\vee}{u_m}\, =\, -\frac{1}{4\pi^2}  h_m 
\frac{ \ \cof \,\big( \A : \sgn \ot \sgn \big)^\top }{ \det \big( \A : \sgn \ot \sgn \big)}\, \overset{\vee}{f},
\]
a.e.\ on $\R^n$. Since $h_m(-z)=h_m(z)$ for all $z\in \R^n$, we get
\begin{align}
\widehat{u_m} (z)\, &=\, -\frac{1}{4\pi^2}  h_m(z)
\frac{ \ \cof\, \Big( \A :\dfrac{-z}{|-z|} \ot \dfrac{-z}{|-z|} \Big)^\top }{ \det \Big(\A : \dfrac{-z}{|-z|}  \ot \dfrac{-z}{|-z|} \Big) }\, \widehat{f}(z) \nonumber\\
&=\, -\frac{1}{4\pi^2}  h_m(z) 
\frac{\ \cof\, \Big( \A : \dfrac{z}{|z|} \ot \dfrac{z}{|z|} \Big)^\top }{ \det \Big(\A : \dfrac{z}{|z|} \ot \dfrac{z}{|z|} \Big) }\, \widehat{f}(z) . \nonumber
\end{align}
Hence, by the identity \eqref{iden}, we deduce
\[
\widehat{u_m} (z)\, =\, -\frac{1}{4\pi^2}  h_m(z) \Big( \A : \dfrac{z}{|z|} \ot \dfrac{z}{|z|} \Big)^{-1}\, \widehat{f}(z) ,
\]
a.e.\ on $\R^n$, which we rewrite as
\beq \label{4.10}
 \A :\widehat{u_m} (z) \ot (2\pi i z) \ot (2\pi i z)\, =\, \big(h_m(z)|z|^2\big)\, \widehat{f}(z).
\eeq
Equivalently,
\beq \label{4.10a}
\A : \widehat{D^2u_m}(z)\, =\, \big(h_m(z)|z|^2\big)\, \widehat{f}(z).
\eeq
By \eqref{4.7a} we have that 
\beq \label{4.14}
0 \, \leq \, h_m(z)|z|^2 \, \leq\, 1
\eeq
and hence by \eqref{4.14}, \eqref{4.10a}, we may employ Proposition \ref{pr2}, Remark \ref{rem} Fourier inversion and Plancherel theorem to infer that each $u_m$ satisfies 
\begin{align}
\| u_m \|_{W^{2,2}_\ast(\Om)}\ &\leq\ C\|D^2u_m\|_{L^2(\Om)} \nonumber\\
&= \ C\Big\| \left[  \big(h_m |\cdot|^2\big) \, \widehat{f} \right]^{\vee} \Big\|_{L^2(\Om)} \nonumber\\
&= \ C\left\|   \big(h_m |\cdot|^2\big)\, \widehat{f} \right\|_{L^2(\Om)} \nonumber\\
&\leq \ C\|  f  \|_{L^2(\Om)}. \nonumber
\end{align}
Hence, $(u_m)_1^\infty$ is bounded in $W^{2,2}_\ast(\R^n)^N$ and as such there is a subsequence of $m$'s and a map $u\in W^{2,2}_\ast(\R^n)^N$ such that, along the subsequence,
\begin{align}
  u_m&\, \lharpoonup u,\ \ \ \text{ in $L^{2^{**}}(\R^n)^N$ as } m\ri \infty \text{ (and a.e.\ on $\R^n$)}, \nonumber\\
 Du_m&\, \lharpoonup Du,\, \text{ in $L^{2^*}(\R^n)^{Nn}$ as } m\ri \infty, \nonumber\\
 D^2u_m&\, \lharpoonup D^2u,\, \text{in $L^{2}(\R^n)^{Nn^2}$ as } m\ri \infty. \nonumber
\end{align}
By \eqref{4.14} and since $h_m(z)|z|^2\ri 1$ for a.e.\ $z\in \R^n$, the Dominated Convergence theorem implies
\[
\big(h_m\, | \cdot |^2 \big)\widehat{f} \larrow \widehat{f}, \ \text{ in  $L^2(\R^n)^N$ as } m\ri \infty. 
\]
By passing to the weak limit as $m\ri \infty$ in \eqref{4.10a}, the Fourier inversion formula implies that the limit $u$ solves 
\[
\A:D^2u\, =\, f
\]
a.e.\ on $\R^n$. By passing to the limit as $m\ri \infty$ in \eqref{4.8}, we obtain the desired representation formula \eqref{4.4a}. Uniqueness of the limit $u$ (and hence independence from the choice of sequence $h_m$) follows from the a priori estimate and linearity. The theorem ensues.       \qed

\section{Existence-uniqueness for fully nonlinear systems} \label{section4}
\ms

We now come to the general fully nonlinear system \eqref{1.1}. We will utilise the results of Sections \ref{section2} and \ref{section3} plus a result of Campanato on near operators, which is recalled later. Our ellipticity condition of Definition \ref{def1} will work as a ``perturbation device", allowing to establish existence for the nonlinear problem by showing it is ``near" a linear well-posed problem. In view of the well-known problems to pass to limits with weak convergence in nonlinear equations, Campanato's idea furnishes an alternative to the stability problem for nonlinear equations, by avoiding this insuperable difficulty. 

The main result of this paper and this section is the next theorem:

\bt[Existence-Uniqueness] \label{th2} Let $n\geq 5$, $N\geq 2$ and let also 
\[
F : \R^n \by \big(\R^{N}\!\ot \mS(n)\big) \larrow \R^N
\]
be a Carath\'eodory map, satisfying Definition \ref{def2} for $\Om=\R^n$ and $F(\cdot,\textbf{0})=0$. Then, for any $f\in L^2(\Om)^N$, the system
\[
F(\cdot,D^2u) \,=\,  f, \ \ \text{ a.e.\ on }\R^n,
\]
has a unique global strong a.e.\ solution $u$ in the space $W^{2,2}_\ast(\R^n)^N$ (given by \eqref{1.10}), which also satisfies the estimate
\beq \label{4.4a}
\|u\|_{L^{2^{**}}(\R^n)}\ +\ \|Du\|_{L^{2^*}(\R^n)}\ +\ \|D^2u\|_{L^{2}(\R^n)}\ \leq \ C\|f\|_{L^2{(\R^n)}},
\eeq
for some $C>0$ depending only on $F$ and the dimensions.  Moreover, for any two maps $w,v \in W^{2,2}_\ast(\R^n)^N$, we have
\beq \label{5.2}
\|w-v\|_{W^{2,2}_\ast(\R^n)}\, \leq\, C \big\|F(\cdot,D^2w)-F(\cdot,D^2v) \big\|_{L^{2}(\R^n)},
\eeq
for some $C>0$ depending only on $F$ and the dimensions. The norm of $W^{2,2}_\ast(\R^n)^N$ is given in Remark \ref{rem}.
\et
 
We note that in view of Lemma \ref{pr1}, the assumption that $F$ satisfies Definition \ref{def2} is equivalent to that $F$ satisfies Definition \ref{def1} plus Lipschitz continuity with respect to the second argument, essentially uniformly with respect to the first argument. We also note that \eqref{5.2} is a strong uniqueness estimate, which is a form of ``comparison principle in integral norms". Moreover, the restriction to homogeneous boundary condition ``$u=0$ at $\infty$" does not harm generality, since the Dirichlet problem we solve is equivalent to a Dirichlet problem with non-homogeneous boundary condition by redefining the nonlinearity $F$ in the standard way.

The proof of Theorem \ref{th2} utilises the following result of Campanato taken from \cite{C5}, whose short proof is given for the sake of completeness at the end of the section:

\bt[Campanato's near operators] \label{th3}  Let $F,A : \mathfrak{X} \larrow X$ be two maps from the set $\mathfrak{X} \neq \emptyset$ to the Banach space $(X,\|\cdot\|)$. Suppose there exists $0<K<1$ such that
\beq \label{5.3}
\Big\| F[u]-F[v]-\big(A[u]-A[v]\big) \Big\| \, \leq\, K \big\| A[u]-A[v] \big\|,
\eeq
for all $u,v \in \mathfrak{X}$. Then, if $A$ is a bijection, $F$ is a bijection as well.
\et

\BPT \ref{th2}. Let $\al$ be the $L^\infty$ function of Definition \ref{def2}. By our assumptions on $F$, Proposition \ref{pr1} implies that there exists $M>0$ depending only on $F$, such that for any $u\in W^{2,2}_\ast(\R^n)^N$, we have
\begin{align} \label{5.5}
\big\|\al(\cdot)F(\cdot,D^2u) \big\|_{L^2(\R^n)} \, &\leq\, \big\|\al(\cdot) F(\cdot,\textbf{0})\big\|_{L^2(\R^n)}\, +\, M\|\al\|_{L^\infty(\R^n)} \| D^2u\|_{L^2(\R^n)} \\
&=\, M\|\al\|_{L^\infty(\R^n)} \| D^2u\|_{L^2(\R^n)} .\nonumber\\
&\leq\, M\|\al\|_{L^\infty(\R^n)} \| u\|_{W^{2,2}_\ast(\R^n)} .\nonumber
\end{align}
The last inequality is a consequence of Remark \ref{rem}. Let also $\A \in \mS(N\!\by\! n)$ be the tensor given by Definition \ref{def2} corresponding to $F$. Then, we have
\beq \label{5.6}
\|\A:D^2u\|_{L^2(\R^n)}\ \leq\ |\A|\, \| D^2u\|_{L^2(\R^n)} \ \leq\,  |\A| \| u\|_{W^{2,2}_\ast(\R^n)} .
\eeq
By \eqref{5.5} and \eqref{5.6} we obtain that the operators 
\beq
\left\{
\begin{array}{l}
A[u]\ :=\ \A :D^2u, \ms\\
F[u]\ :=\ \al(\cdot)F(\cdot, D^2u),
\end{array}
\right.
\eeq
map $W^{2,2}_\ast(\R^n)^N$ into  $L^2(\R^n)^N$. Let $u,v \in W^{2,2}_\ast(\R^n)^N$. By Definition \ref{def2} and the a priori hessian estimate of Proposition \ref{pr2} we have
\begin{align} \label{5.7}
\Big\|\al(\cdot)\Big(F(\cdot, & D^2u) -  F(\cdot,D^2v)\Big)  -\A: \big(D^2u-D^2v \big)\Big\|^2_{L^2(\R^n)} \nonumber\\
& \leq\ \be  \big\|\A:(D^2u-D^2v)\big\|^2_{L^2(\R^n)}\ + \ \ga \big\|\A:(D^2u-D^2v)\big\|^2_{L^2(\R^n)}  \\
& \leq\ (\be +\ga) \big\|\A:(D^2u-D^2v)\big\|^2_{L^2(\R^n)}. \nonumber
\end{align}

Theorem \ref{th1} implies that the linear operator
\[
A\ :\  W^{2,2}_\ast(\R^n)^N \larrow  L^2(\R^n)^N 
\]
is a bijection. Hence, in view of the inequality \eqref{5.7} and the fact that $\sqrt{\be+\ga}<1$, Campanato's Theorem \ref{th3} implies that $F:  W^{2,2}_\ast(\R^n)^N \larrow  L^2(\R^n)^N $ is a bijection as well. As a result, for any $g\in L^2(\R^n)^N$, the PDE system $\al(\cdot)F(\cdot,D^2u)=g$ has a unique solution in $W^{2,2}_\ast(\R^n)^N$. Since $\al,1/\al \in L^\infty(\R^n)$, by selecting $g=\al f$, we conclude that the problem \eqref{1.1} has a unique solution  in $W^{2,2}_\ast(\R^n)^N$. Finally, by \eqref{5.7} we have
\[
\Big\|F(\cdot,D^2u) -  F(\cdot,D^2v)\Big\|_{L^2(\R^n)} \,  \geq \ \frac{1-\sqrt{\be+\ga}}{\|\al\|_{L^\infty(\R^n)}} \, \big\| \A:(D^2u -D^2v) \big\|_{L^2(\R^n)}
\]
and by Proposition \ref{pr2} and Remark \ref{rem}, we deduce the estimate
\begin{align} 
\Big\|F(\cdot,D^2u) -  F(\cdot,D^2v)\Big\|_{L^2(\R^n)} \, &\geq \ \left(\nu(\A)\frac{1-\sqrt{\be+\ga}}{\|\al\|_{L^\infty(\R^n)}}\right) \big\|D^2u -D^2v\big\|_{L^2(\R^n)} \nonumber\\
& \geq\ C\, \| u -v \|_{W^{2,2}_\ast(\R^n)},\nonumber
\end{align}
for some $C>0$. The theorem ensues.      \qed

\ms

We conclude this section with the proof of Campanato's theorem on near operators taken from \cite{C5}, which we provide for the convenience of the reader. 

\BPT \ref{th3}. It suffices to show that for any $f\in X$, there is a unique $u\in \mathfrak{X}$ such that
\[
F[u]\, =\, f.
\]
In order to prove that, we first turn $\mathfrak{X}$ into a complete metric space, by pulling back the structure from $X$ via $A$: for, we define the distance
\[
d(u,v)\, :=\, \big\| A[u]-A[v]\big\|.
\]
Next, we fix an $f\in X$ and define the map
\[
T\ : \ \mathfrak{X} \larrow \mathfrak{X}\ , \ \ \ T[u]\, :=\, A^{-1}\Big(A[u]-\big(F[u]-f \big) \Big).
\]
We conclude by showing that $T$ is a contraction on $(\mathfrak{X},d)$, and hence has a unique $u\in \mathfrak{X}$ such that $T[u]=u$. The latter equality is equivalent to $F[u]=f$, and then we will be done.  Indeed, we have that
\begin{align}
d\Big( T[u],T[v] \Big)  \, &=\ \Big\|  \left(A[u]-\big(F[u]-f \big) \right) \, -\,  \left(A[v]-\big(F[v]-f \big) \right) \Big\| \nonumber \\
 &=\ \Big\| A[u] -A[v] -\big(F[u]  - F[v]\big) \Big\|, \nonumber 
\end{align}
and hence
\begin{align}
d\Big( T[u],T[v] \Big)  \ \, 
&\!\!\!\!\overset{\eqref{5.3}}{\leq}   K \big\| A[u]-A[u] \big\| \nonumber\\
  &= \  K \, d(u,v). \nonumber
\end{align}
Since $K<1$, the conclusion follows and the theorem ensues.              \qed

\section{Extensions} \label{section5}

In this section we discuss an extension of Theorem \ref{th2} in the form of ``stability theorem for strong solutions".

\bt[Stability of strong solutions] \label{th3} Let $n\geq 5$, $N\geq 2$ and $F,G : \R^n \by \big(\R^{N}\!\ot \mS(n)\big) \larrow \R^N$ Carath\'eodory maps. We suppose that 
\[
F\ :\  W^{2,2}_\ast(\R^n)^N \larrow  L^2(\R^n)^N 
\]
is a bijection, where the space $W^{2,2}_\ast(\R^n)^N$ is given by \eqref{1.10}.  If $G(\cdot,\textbf{0})=0$ and
\beq \label{6.1}
\underset{x\in \R^n}{\ess\,\sup}\sup_{\X\neq\Y}\left| \frac{\big(F(x,\Y)-F(x,\X)\big)-\big(G(x,\Y)-G(x,\X) \big)}{|\Y \, -\, \X|}\right| \ < \ \nu(F)
\eeq
where
\beq \label{6.2}
\nu(F)\ :=\ \inf_{v\neq w}\frac{\big\|F(\cdot,D^2w)-F(\cdot,D^2v) \big\|_{L^{2}(\R^n)}}{\|D^2w-D^2v\|_{L^2(\R^n)}}\ >\ 0,
\eeq
then, for any given $g\in L^2(\R^n)^N$, the system
\[
G(\cdot,D^2u) \,=\,  g, \ \ \text{ a.e.\ on }\R^n,
\]
has a unique global strong a.e.\ solution $u$ in the space $W^{2,2}_\ast(\R^n)^N$.
\et
 
Theorem \ref{th2} provides sufficient conditions on $F$ is order to obtain solvability. Hence, every $G$ which is ``close to $F$" in the sense of \eqref{6.1}, gives rise to a nonlinear coefficient such that the respective global Dirichlet problem is uniquely solvable.

\BPT \ref{th3}. 
We denote the right hand side of \eqref{6.1} by $\nu(F,G)$ and we may rewrite  \eqref{6.1}  as
\beq \label{6.3}
0\,<\,  \nu(F,G)\, <\, \nu(F).
\eeq
For any $u,v\in W^{2,2}_\ast(\R^n)^N$, we have
\begin{align} 
\Big\| F& (\cdot,D^2u) -  F(\cdot,D^2v)  -\big( G(\cdot,D^2u) -  G(\cdot,D^2v)\big) \Big\|_{L^2(\R^n)} \nonumber\\
& \leq \left( \underset{\R^n}{\ess\,\sup} 
\sup_{\X\neq \Y}\,
\left| \frac{ 
F(\cdot,\Y) -F(\cdot,\X) -
\big( G(\cdot,\Y) -G(\cdot,\X) \big) }{|\Y \, -\, \X|}
\right| \right)
\big\| D^2u-D^2v\big\|_{L^2(\R^n)}  \nonumber\\
& = \, \nu(F,G) \big\| D^2u-D^2v\big\|_{L^2(\R^n)} \nonumber\\
&\leq\, \frac{\nu(F,G)}{\nu(F)} \big\|F(\cdot,D^2u) -F(\cdot,D^2v)\big\|_{L^2(\R^n)}. \nonumber
\end{align}
Hence, we obtain the inequality
\begin{align}  \label{6.4}
\Big\| F (\cdot,D^2u) -  F(\cdot,D^2v)  -&\big( G(\cdot,D^2u) -  G(\cdot,D^2v)\big) \Big\|_{L^2(\R^n)} \nonumber\\
&\leq\, \frac{\nu(F,G)}{\nu(F)} \big\|F(\cdot,D^2u) -F(\cdot,D^2v)\big\|_{L^2(\R^n)},
\end{align}
which is valid for any $u,v\in W^{2,2}_\ast(\R^n)^N$. By \eqref{6.3}, Remark \ref{rem} and the inequality above for $v\equiv 0$, we have that $F,G$ map $W^{2,2}_\ast(\R^n)^N$ into  $L^2(\R^n)^N$. By assumption, $F :  W^{2,2}_\ast(\R^n)^N \larrow  L^2(\R^n)^N$ is a bijection. Hence, in view of Campanato's Theorem \ref{th3}, the inequality \eqref{6.4} implies that $G :  W^{2,2}_\ast(\R^n)^N \larrow  L^2(\R^n)^N$ is a bijection as well.  The theorem ensues.     \qed

\section{Motivations and Potential Applications} \label{Motivations_and_Applications}

In this section we collect some material relevant to the problem we are considering in this paper  and to which our results may apply by perhaps imposing appropriate conditions and/or restrictions. Our motivation to study this problem comes from the necessity to understand PDE systems arising in Differential Geometry, Mathematical Physics and Calculus of Variations: 

\ms

\noi \textbf{The Harmonic map problem:} Given two Riemannian manifolds $(\mathcal{M},\gamma)$, $(\mathcal{N},g)$, then a smooth map $u:\mathcal{M}\larrow \mathcal{N}$ is called a harmonic map if and only if the following PDE system is satisfied:
\[
\Delta_{M} u_\al \,+\, \gamma^{ij}\,\Gamma^{\al}_{\mu \nu}(u)\, D_i u_\nu D_j u_\mu \, =\, 0. 
\]
Here $\Gamma$ denote the Christoffel symbols of the target metric $g$ and 
\[
\Delta_{\mathcal{M}} u\, =\, \frac{1}{det(\gamma)}\Div \left(\sqrt{\det(\gamma)}\gamma\, D u\right)
\]
is the Laplace-Beltrami operator.  If $u$ additionally is an isometric embedding, then $u(\mathcal{M})$ is a minimal submanifold of $\mathcal{N}$. The problem is highly non-trivial even when $\mathcal{N}=\R^N$ and in the case of codimension greater than one, that is when \[
\dim(\mathcal{N})\, -\, \dim(\mathcal{M})\, \geq \, 2
\]
it is far from well understood. Moreover, it is well known (see e.g. Lawson-Osserman \cite{LO}) that then (by using the properties of the second fundamental form) the system above can be written in an equivalent formulation of a non-divergence 2nd order elliptic system of the form we are considering in this paper. 

\ms

\noi \textbf{Elliptic problem involving the Ricci curvature:} Let $(\mathcal{M},g)$ be a Riemannian metric. Then the principal part ``P.P." of the Ricci curvature in coordinates is
\[
P.P.(R_{\mu\nu})= g^{\alpha\beta} \Big(D_{\mu}D_{\alpha}g_{\beta\nu}+D_{\nu}D_{\alpha}g_{\beta\mu}-D_{\mu}D_{\nu}g_{\alpha\beta}-D_{\alpha}D_{\beta}g_{\mu\nu} \Big) .
\]
If we now consider harmonic coordinates $(x_1,...,x_n)$, that is those for which $\Delta_\mathcal{M} x_i =0$, then it follows from standard computations in Riemannian geometry that the expression
\[
D_{\mu}D_{\alpha}g_{\beta\nu}+D_{\nu}D_{\alpha}g_{\beta\mu}-D_{\mu}D_{\nu}g_{\alpha\beta}
\] 
is given by an expression which involves at most one derivative of the metric components. Therefore, for harmonic coordinates we have
\[
P.P.(R_{\mu\nu})\, =\, g^{\alpha\beta}D_{\alpha}D_{\beta}g_{\mu\nu}  
\]
and hence, any identity that the Ricci curvature satisfies can be easily seen to correspond to a non-divergence 2nd order elliptic system for the metric components. 

\ms

\noi \textbf{Elliptic problems arising in the Einstein equations and in Conformal Geometry:} The celebrated equations in the vacuum in local coordinates read
\[
R_{\al \be}\, =\, 0,
\]
If the (unknown) metric $g$ admits a Killing vector field which is timelike, then the Einstein equations reduce by using the above to a quasilinear non-divergence elliptic system (see e.g.\ \cite{HE}). Moreover, it is also well known that fully nonlinear elliptic elliptic systems of the type we consider herein arise in Conformal Geometry, see e.g. Trudinger \cite{Tr}, \cite{G}.

\ms

\noi \textbf{Non-convex 2nd order Variational problems:} Consider the functional
\[
E(u,\Om)\, =\, \int_\Om f\big( x,D^2u(x)\big)\, dx
\] 
placed in the space $W^{2,p}(\Om)^N$. Then, it is well known in Calculus of Variations (see e.g.\ \cite{D, DM, GM}) that if $f$ fails to be quasiconvex with respect to the Hessian argument, then minimisers may well not exist in the Sobolev space. Then, a standing idea in order to construct minimisers is to solve a fully nonlinear 2nd order PDE with vectorial solution of the form 
\[
f(\cdot,D^2u)\, =\, h
\]
on the subdomain of $\Om$ obtained when we consider the set whereon $f$ is strictly greater than its quasiconvex envelope $\overline{f}$.

\ms

\noi \textbf{The equations of vectorial $L^\infty$ variational problems:}  Calculus of Variations in $L^\infty$ has a long history and was pioneered by Aronsson in the 1960s (see \cite{A1}-\cite{A7}) who was the first to consider variational problems for supremal functionals of the form
\[
E_\infty(u,\Om)\, =\, \big\| H(\cdot,u,Du) \big\|_{L^\infty(\Om)}.
\]
However, until the early 2010s the field was essentially restricted to the scalar case. The foundations of the vectorial case have been laid in a series of recent papers of the author (see \cite{K1}-\cite{K7}). In the simplest case of 
\[
E_\infty(u,\Om)\, =\, \| Du \|_{L^\infty(\Om)}
\]
applied to Lipschitz maps $u : \Om \sub \R^n \larrow \R^N$, the counterpart of the ``Euler-Lagrange" equations is the so-called $\infty$-Laplace system:
\[
\De_\infty u \, :=\, \Big(Du \ot Du\, +\, |Du|^2[Du]^\bot \! \ot I \Big):D^2u\, =\, 0.
\]
In the above, $[Du(x)]^\bot$ denotes the orthogonal projection on the nullspace of the operator $Du(x)^\top : \R^N \larrow \R^n$ and in index form reads
\[
\Big(D_iu_\al \, D_ju_\be\, +\, |Du|^2[Du]_{\al \be}^\bot \, \de_{ij}\Big)\, D_{ij}^2u_\be\, =\, 0.
\]
The system above is nondivergence quasilinear degenerate elliptic, has discontinuous coefficients and behaves like a fully nonlinear elliptic system. The problem we consider herein with pure Hessian dependence is an essential stepping stone for the understanding on the $\infty$-Laplace system. Indeed, the results of this paper have been invaluable tools in the very recent papers \cite{K10, K11} wherein we study the $L^\infty$ equations.

\ms

\ms

\noi \textbf{Acknowledgement.} I would like to thank Tristan Pryer for our inspiring scientific discussions. I am also grateful to Stefanos Aretakis for our discussions on elliptic problems arising in Differential Geometry.

\ms

\ms

\bibliographystyle{amsplain}

\end{document}